\documentclass[11pt, twoside]{article}

\usepackage{extarrows}
\usepackage{enumerate}
\usepackage{amsmath}
\usepackage{amssymb}
\usepackage{amsthm}
\usepackage[all]{xypic}
\usepackage[all,cmtip,poly]{xy}
\usepackage{latexsym}
\usepackage{pstricks,hyperref}
\usepackage{amsbsy}
\usepackage{amscd}
\usepackage{mathrsfs}
\usepackage{txfonts}
\usepackage{amsfonts}
\usepackage{graphicx}
\usepackage[top=1.40in,bottom=1.30in,left=1.20in,right=1.20in]{geometry}
\newtheorem{thm}{Theorem}[section]
\newtheorem{cor}[thm]{Corollary}
\newtheorem{lem}[thm]{Lemma}

\newtheorem{defn}[thm]{Definition}
\newtheorem{rem}[thm]{Remark}

\usepackage{pstricks,color}

\usepackage{indentfirst}
\begin{document}

\begin{center}
{\Large \bf A note on abelian quotient categories\footnotetext{This work was supported by the Hunan Provincial Natural Science Foundation of China (Grants No. 2018JJ3205) and the NSF of China (Grants No. 11671221)}}

\bigskip

{\large Panyue Zhou}
\bigskip
\end{center}

\def\s{\stackrel}
\def\Longrightarrow{{\longrightarrow}}
\def\A{\mathcal{A}}
\def\B{\mathcal{B}}
\def\C{\mathscr{C}}
\def\D{\mathsf{D}}
\def\T{\mathcal{T}}
\def\M{\mathcal{M}}
\def\E{\mathcal{E}}
\def\R{\mathcal{R}}
\def\S{\mathcal{S}}
\def\H{\mathcal{H}}
\def\U{\mathscr{U}}
\def\V{\mathscr{V}}
\def\W{\mathscr{W}}
\def\X{\mathscr{X}}
\def\Y{\mathscr{Y}}
\def\Z{\mathcal {Z}}
\def\I{\mathcal {I}}
\def\RR{\mathcal{R}\ast\mathcal{R}[1]}
\def\Aut{\mbox{Aut}}
\def\coker{\mbox{coker}}
\def\Ker{\mbox{Ker}}
\def\deg{\mbox{deg}}
\def\dim{\mbox{dim}}
\def\Ext{\mbox{Ext}}
\def\Hom{\mbox{Hom}}
\def\DHom{\mbox{DHom}}
\def\Gr{\mbox{Gr}}
\def\id{\mbox{id}}
\def\Im{\mbox{Im}}
\def\ind{\mbox{ind}}
\def\Int{\mbox{Int}}
\def\ggz{\Gamma}
\def\la{\Lambda}
\def\bz{\beta}
\def\az{\alpha}
\def\gz{\gamma}
\def\da{\delta}
\def\fs{{\mathfrak{S}}}
\def\ff{{\mathfrak{F}}}
\def\zz{\zeta}
\def\thz{\theta}
\def\raw{\rightarrow}
\def\ole{\overline}
\def\cat{\C_{F^m}(\H)}
\def\fun{F_\la}
\def\sttm{\mbox{s}\tau\mbox{-tilt}\la}
\def\wte{\widetilde}
\def \text{\mbox}
\def\Mod{\mathsf{Mod}}
\hyphenation{ap-pro-xi-ma-tion}
\newcommand{\pd}{\mathsf{pd}\hspace{.01in}}
\newcommand{\add}{\mathsf{add}\hspace{.01in}}
\newcommand{\Fac}{\mathsf{Fac}\hspace{.01in}}
\newcommand{\thick}{\mathsf{thick}\hspace{.01in}}
\newcommand{\End}{\operatorname{End}\nolimits}
\renewcommand{\Mod}{\mathsf{Mod}\hspace{.01in}}
\renewcommand{\mod}{\mathsf{mod}\hspace{.01in}}
\newcommand{\proj}{\mathsf{proj}\hspace{.01in}}
\newcommand{\co}{\mathsf{Coker}\hspace{.01in}}
\newcommand{\pr}{\mathsf{pr}\hspace{.01in}}
\begin{abstract}
 Let $\C$ be a triangulated category with a Serre functor $\mathbb{S}$ and $\X$ a
non-zero contravariantly finite rigid subcategory of $\C$. Then $\X$ is cluster tilting
if and only if the quotient category $\C/\X$ is abelian and $\mathbb{S}(\X)=\X[2]$.
As an application, this
result generalizes work by Beligiannis.\\[0.2cm]
\end{abstract}

\medskip

\pagestyle{myheadings}
\markboth{\rightline {\scriptsize   P. Zhou}}
         {\leftline{\scriptsize A note on abelian quotient categories}}
\thispagestyle{plain}

\section{Introduction}
Let $\C$ be an additive category and $\X$ be a subcategory of $\C$ and let $M\in \C$. A morphism $f_M\colon X_M \to M$ is called a right $\X$-approximation of $M$ if $X_M\in\X$ and every morphism from an object in $\X$ to $M$
to factors through $f_M$. The subcategory $\X$ is said to be \emph{ contravariantly finite} in $\C$, if every object in $\C$ has a right $\X$-approximation.  A \emph{left $\X$-approximation} and a \emph{covariantly finite subcategory} of $\C$
are dually defined. A contravariantly and covariantly finite subcategory is called
\emph{functorially finite}.   For more details, see \cite{AR}.

Recall the notion of cluster tilting subcategories from \cite[Definition 3.1]{KZ}.

\begin{defn}\emph{\cite[Definition 3.1]{KZ}}
Let $\C$ be a triangulated category with a shift functor $[1]$.
\begin{enumerate}
\item[\emph{(1)}] A subcategory $\X$ of $\C$ is called {\rm rigid} if ${\rm Hom}_{\C}(\X, \X[1])=0$.
\item[\emph{(2)}] A functorially finite subcategory $\X$ of $\C$ is called {\rm cluster tilting} if
$$\X=\{M\in\C\ |\ {\rm Hom}_{\C}(\X, M[1])=0\}=\{M\in\C\ |\ {\rm Hom}_{\C}(M, \X[1])=0\}.$$
\end{enumerate}
\end{defn}

\begin{rem}\label{rem:ctsubcat}
In fact, Koenig and Zhu \emph{\cite[Lemma 3.2]{KZ}} indicate that a subcategory $\X$ of a triangulated categoru $\C$ is cluster tilting if and only if it is contravariantly finite in $\C$ and
$$\X=\{M\in\C\ |\ {\rm Hom}_{\C}(\X, M[1])=0\}.$$
\end{rem}

Koenig and Zhu gave a general framework from passing from triangulated categories
to abelian categories by factoring out cluster tilting subcategories. More precisely, they proved the following.

\begin{thm}\label{thm1}\emph{\cite[Theorem 3.3]{KZ}}
Let $\C$ be a triangulated category and $\X$ a
cluster tilting subcategory of $\C$.
Then the quotient category $\C/\X$ is an abelian category.
\end{thm}

A triangulated category is called \emph{connected} if it cannot be decomposed into a direct sum of two non-zero triangulated subcategories. Beligiannis showed the following
characterization of cluster tilting subcategories which complements, and was inspired by Koenig and Zhu.

\begin{thm}\label{thm2}\emph{\cite[Theorem 7.3]{Be}}
 Let $\C$ be a \textbf{connected} triangulated category with a Serre functor $\mathbb{S}$ and $\X$ a
non-zero \textbf{functorially finite} rigid subcategory of $\C$. Then $\X$ is cluster tilting
if and only if the quotient category $\C/\X$ is abelian and $\mathbb{S}(\X)=\X[2]$.
\end{thm}

Our aim in this article is to give a simple short proof of the following more general
result.

\begin{thm}
 Let $\C$ be a triangulated category with a Serre functor $\mathbb{S}$ and $\X$ a
non-zero \textbf{contravariantly finite} rigid subcategory of $\C$. Then $\X$ is cluster tilting
if and only if the quotient category $\C/\X$ is abelian and $\mathbb{S}(\X)=\X[2]$.
\end{thm}

\section{Preliminaries}
We recall the definition of Auslander-Reiten triangle.
\begin{defn}{\emph{\cite[Definition 4.1]{Ha}}}
Let $\C$ be a triangulated category. A triangle
$$A\xrightarrow{~u~}B\xrightarrow{~v~}C\xrightarrow{~w~}A[1]$$
in $\C$ is called an \emph{Auslander-Reiten triangle} if
\begin{itemize}
\item $A,B$ are indecomposable;

\item $w\neq 0$

\item If $f\colon X\to C$ is not a split epimorphism, then there exists a morphism
$f'\colon X\to B$ such that $f'v=f$.
\end{itemize}
In this case, we call $u$ is a source morphism and $v$ is  a sink morphism.  The notions of a source morphism and a sink morphism are also known as minimal left almost split morphism and minimal right split morphism, respectively.
\end{defn}

 Assume that $\C$ is a $k$-linear Hom-fnite triangulatedcategory where $k$ is a field. Recall from \cite{BK} a Serre functor $\mathbb{S}\colon \C\to \C$ is a $k$-linear equivalence with bifunctorial isomorphisms
$${\rm Hom}_{\C}(A, B)\simeq D{\rm Hom}_{\C}(B, \mathbb{S}A)$$
for any $A, B\in \C$, where $D$ is the duality over $k$. Reiten and Van den Bergh \cite{RV} proved that if $\C$
admits a Serre functor $\mathbb{S}$, then $\C$ has Auslander-Reiten triangles. Moreover, if $\tau$ is the Auslander-Reiten translation in $\C$, then $\mathbb{S}\simeq\tau[1]$.
 We say that a triangulated category $\C$ is \emph{$2$-Calabi-Yau} if $\mathbb{S}\simeq [2]$.

Let $\C$ be an additive category and $\X$ be a subcategory of $\C$.
We denote by $\C/\X$
the category whose objects are objects of $\C$ and whose morphisms are elements of
$\Hom_{\C}(A,B)/\X(A,B)$ for $A,B\in\C$, where $\X(A,B)$ the subgroup of $\Hom_{\C}(A,B)$ consisting of morphisms
which factor through an object in $\X$.
Such category is called the \emph{quotient category}
of $\C$ by $\X$. For any morphism $f\colon A\to B$ in $\C$, we denote by $\overline{f}$ the image of $f$ under
the natural quotient functor $\C\to\C/\X$.

\section{Main result}
In order to prove our main result, we need the following some lemmas.

\begin{lem}\label{y1}\emph{\cite[Theorem 3.1]{BM}} and \emph{\cite[Theorem 2.2]{J}} 
 Let $\C$ be a triangulated category and $\X$ a contravariantly finite
 subcategory of $\C$. Then the quotient category $\C/\X$ is a left triangulated category with respect to the following loop functor and left triangles:
 \begin{itemize}
 \item[\emph{(a)}] For any object $C\in\C$, we take a triangle
 $$\Omega C\xrightarrow{~\alpha~}X\xrightarrow{~\beta~}C\xrightarrow{~\gamma~}\Omega C[1]$$
 where $\beta$ is a right $\X$-approximation of $C$. Then $\Omega$ gives a well-defined functor of $\C/\X$, which is the loop functor of $\C/\X$.
 \item[\emph{(b)}] A left triangle in $\C/\X$ is by definition a diagram which is isomorphic in $\C/\X$ to a
diagram
$$\Omega C\xrightarrow{~\overline{u}~}A\xrightarrow{~\overline{v}~}B\xrightarrow{~\overline{w}~}C$$
arising by forming the homotopy pull-back in $\C$
\begin{equation}\label{tb}
\begin{array}{l}
\xymatrix{
\Omega C \ar[r]^{u}\ar@{=}[d]&A \ar[r]^{v} \ar[d] &B \ar[d]^{w}\ar[r]&\Omega C[1]\ar@{=}[d] \\
\Omega C \ar[r]^{\alpha} &X \ar[r]^{\beta} &C \ar[r]^{\gamma}&\Omega C[1]}
\end{array}
\end{equation}
of a morphism $w\colon B\to C$ along the triangle $\Omega C\xrightarrow{~~}X\xrightarrow{~g~}C\xrightarrow{~~}\Omega C[1]$ where $\beta$ is a right $\X$-approximation of $C$.
Equivalently it is easy to see that the left triangles in $\C/\X$ are the
diagrams which are isomorphic in $\C/\X$ to diagrams \emph{(\ref{tb})}, arising from triangles
$A\xrightarrow{~~}B\xrightarrow{~~}C\xrightarrow{~~}A[1]$ in $\C$, where the morphism 
$\emph{\Hom}_{\C}(\X,B)\to \emph{\Hom}_{\C}(\X,C)$ is an epimorphism.

 \end{itemize}
\end{lem}

\begin{lem}\label{y2}
 Let $\C$ be a triangulated category and $\X$ a
 subcategory of $\C$. If
$$A\xrightarrow{~f~}B\xrightarrow{~g~}C\xrightarrow{~h~}A[1]$$
is the Auslander-Reiten triangle in $\C$ and $A\notin\X$, then we have the following
exact sequence in the quotient category $\overline{\C}:=\C/\X$:
$$\emph{\Hom}_{\overline{\C}}(C,M)\xrightarrow{~\overline{g}~}\emph{\Hom}_{\overline{\C}}(B,M)\xrightarrow{~\overline{f}~}\emph{\Hom}_{\overline{\C}}(A,M),$$
where $M\in\overline{\C}$.
\end{lem}

\proof Since $gf=0$ and then $\overline{f}\circ \overline{g}=0$. Thus we have ${\rm Im}(\overline{g})\subseteq {\rm Ker}(\overline{f})$.
 \smallskip

 Conversely, let $\overline{u}\colon B\to M$ be any morphism in $\overline{\C}$ such that $\overline{u}\circ\overline{f}=0$.
Then there exists morphisms $s\colon A\to X$ and $t\colon X\to M$ such that $uf=ts$.
Since $A\in\notin\X$, we have that $s\colon A\to X$ is not split monomorphism. So there exist a morphism
$x\colon B\to X$ such that $s=xf$ and then $$(u-tx)f=uf-ts=0.$$
Thus there exists a morphism $y\colon C\to M$ such that $u-tx=yg$.
It follows that $\overline{u}=\overline{y}\circ\overline{g}$.

This shows that ${\rm Ker}(\overline{f})\subseteq{\rm Im}(\overline{g})$.  \qed
\medskip

Now we prove the main result of this article. The proof is not too far from Beligiannis case, compare with \cite[Theorem 7.3]{Be}.

\begin{thm}
 Let $\C$ be a triangulated category with a Serre functor $\mathbb{S}$ and $\X$ a
non-zero contravariantly finite rigid subcategory of $\C$. Then $\X$ is cluster tilting
if and only if the quotient category $\C/\X$ is abelian and $\mathbb{S}(\X)=\X[2]$.
\end{thm}

\proof   We first show the `only if' part.

By Theorem \ref{thm1}, we know that the quotient category $\C/\X$ is an abelian category.

Now we show that $\mathbb{S}(\X)=\X[2]$.

Since $\Hom_{\C}(\X,\mathbb{S}\X[-2][1])=\Hom_{\C}(\X,\mathbb{S}\X[-1])\simeq D\Hom_{\C}(\X[-1],\X)=0$ and $\X$ is cluster tilting,  we have $\mathbb{S}(\X)\subseteq\X[2]$.

Conversely, note that $\Hom_{\C}(\mathbb{S}^{-1}\X[2],\X[1])\simeq D\Hom_{\C}(\X[1],\X[2])=0$, thus
for any $M\in \mathbb{S}^{-1}\X[2]$, we have $\Hom_{\C}(M,\X[1])=0$.

Since $\X$ is contravariantly finite in $\C$, we can take a triangle
\begin{equation}\label{t0}
\begin{array}{l}
N\xrightarrow{~f~}X_0\xrightarrow{~g~}M\xrightarrow{~h~}N[1]
\end{array}
\end{equation}
with $g$ is a right $\X$-approximation of $M$.  Apply the functor $\Hom_{\C}(\X,-)$ to the above triangle, we have
the following exact sequence:
$$\Hom_{\C}(\X,X_0)\xrightarrow{~\Hom_{\C}(\X,g)~}\Hom_{\C}(\X,M)\xrightarrow{~~}\Hom_{\C}(\X,N[1])\xrightarrow{~~}\Hom_{\C}(\X,X_0[1])=0.$$
Since $g$ is a right $\X$-approximation of $M$, we have that $\Hom_{\C}(\X,g)$ is an epimorphism. It follows that
$\Hom_{\C}(\X,N[1])=0$ implies $N\in\X$ since $\X$ is cluster tilting.
We obtain that $h=0$. Thus the triangle (\ref{t0}) splits. It follows that $M$ is a direct summand of $X_0$ and then $M\in\X$.
Hence $M\in \mathbb{S}^{-1}\X[2]\subseteq\X$ and then $\X[2]\subseteq\mathbb{S}\X$.
\medskip

To prove the `if' part.

Since $\X$ is rigid, we have that $\X\subseteq\{M\in\C~|~\Hom_{\C}(\X,M[1])=0\}$.

Suppose that $M$ is an indecomposable object satisfying $\Hom_{\C}(\X,M[1])=0$.

Now we assume $M\notin \X$. Since $$\Hom_{\C}(\X,M[1])\simeq D\Hom_{\C}(M[1],\mathbb{S}\X)\simeq D\Hom_{\C}(M,\mathbb{S}\X[-1])=D\Hom_{\C}(M,\X[1]),$$
 we have $\Hom_{\C}(M,\X[1])=0$. Then $\tau M\notin\X$ since $\Hom_{\C}(M,\tau M[1])\simeq D\Hom_{\C}(M,M)\neq 0$.
 Let
$$\tau M\xrightarrow{~f~}C\xrightarrow{~g~}M\xrightarrow{~h~}\tau M[1]$$
be the Auslander-Reiten triangle ending at $M$.

Since $\X$ is contravariantly finite in $\C$, then for the object $M$ there exists a triangle
$$\Omega M\xrightarrow{~x~}X_0\xrightarrow{~y~}M\xrightarrow{~z~}\Omega M[1]$$
where $y$ is a right $\X$-approximation of $M$. Since $M\notin\X$, we have that $y$ is not a split epimorphism.
It follows that $y$ factors through $g$ and there exists a morphism of triangles
$$\xymatrix{
\Omega M \ar[r]^{x}\ar[d]^{\phi}&X_0 \ar[r]^{y} \ar[d]^{\varphi} &M \ar@{=}[d]\ar[r]^{z\quad}&\Omega M[1]\ar[d] \\
\tau M \ar[r]^{f} &C \ar[r]^{g} &M \ar[r]^{h}&\tau M}
$$
which by Lemma \ref{y1} induces the following left triangle in the left triangulated category $\C/\X$:
$$\Omega M\xrightarrow{~\overline{\phi}~}\tau M\xrightarrow{~\overline{f}~}C\xrightarrow{~\overline{g}~}M$$

If $\overline{g}=0$, then the above left triangle gives a direct sum decomposition
$\tau M\simeq \Omega M\oplus C$
in $\C/\X$. Since $M$ is indecomposable, we have that $\tau M$ is indecomposable.
Thus we have $C=0$ in $\C/\X$ or else $\Omega M=0$ in $\C/\X$, that is to say,
either $C\in\X$ or else $\Omega M\in\X$.
If $\Omega M\in\X$, we have $z\in\Hom_{\C}(M,\X[1])=0$. Then $y$ is a split epimorphism,
namely $M$ lies in $\X$ as a direct summand of $X_0$ and this is a contradiction to $M\notin\X$.
We infer that $C\in\X$.

Since $\overline{g}=0$, there exist morphisms $s\colon C\to X_1$ and $t\colon X_1\to M$ such that
$g=ts$ where $X_1\in\X$. Since $M\notin\X$, we have that $t$ is not a split epimorphism. So there exists
a morphism $k\colon X_1\to C$ such that $t=gk$ and then $g(1-ks)=g-gks=g-ts=0$. Thus there exists a morphism
$m\colon C\to \tau M$ such that $1-ks=fm$ and then $\overline{f}\circ\overline{m}=\overline{1}$. 
We obtain $\overline{f}$ is a split epimorphism in $\C/\X$.
Then $C$ is a direct summand of $\tau M$. Since $\tau M\notin\X$, we have $C\notin\X$. This is a contradiction to
$C\in\X$.

So we get $\overline{g}\neq 0$. Let
\begin{equation}\label{t2}
\begin{array}{l}
M\xrightarrow{~u=\tau^{-1}f~}\tau^{-1}C\xrightarrow{~v=\tau^{-1}g~}\tau^{-1}M\xrightarrow{~w=\tau^{-1}h~}M[1]
\end{array}
\end{equation}
be the Auslander-Reiten triangle starting at $M$. Then $\overline{v}\neq 0$ since $\overline{g}\neq 0$.

We claim that $\overline{v}$ is an epimorphism $\C/\X$. In fact, let $\overline{\alpha}\colon \tau^{-1}M\to L$
be any morphism in $\C/\X$ such that $\overline{\alpha}\circ\overline{v}=0$.
Since $\C/\X$ is abelian, we assume that $\overline{\beta}\colon \tau^{-1}C\to L$ is the cokernel of $\overline{u}$.
Then there exists a morphism $\overline{\omega}\colon L\to \tau^{-1}M$ such that $\overline{v}=\overline{\omega}\circ\overline{\beta}$ since $\overline{v}\circ\overline{u}=0$.
Note that $\overline{\beta} \circ\overline{u}=0$, by Lemma \ref{y2}, there exists a morphism
$\overline{\delta}\colon \tau^{-1}M\to L$ such that $\overline{\beta}=\overline{\delta}\circ\overline{v}$.
It follows that $$\overline{\beta}=\overline{\delta}\circ\overline{v}=\overline{\delta}\circ\overline{\omega}\circ\overline{\beta}
~~\textrm{and }~\overline{v}=\overline{\omega}\circ\overline{\beta}=\overline{\omega}\circ\overline{\delta}\circ\overline{v}.$$
Since $\overline{\beta}$ is an epimorphism in $\C/\X$, we have $\overline{\delta}\circ\overline{\omega}=\overline{1}$
and $\overline{\omega}\circ\overline{\delta}$ is an idempotent in ${\rm End}_{\C/\X}(\tau^{-1}M)$.
Since $\overline{v}\neq 0$, we have $\overline{\omega}\circ\overline{\delta}\neq 0$.
Thus $\overline{\omega}\circ\overline{\delta}$ is not nilpotent.
Since ${\rm End}_{\C/\X}(\tau^{-1}M)$ is local, we have $\overline{\omega}\circ\overline{\delta}$ is an isomorphism.
Hence $\overline{\omega}$ is an isomorphism in $\C/\X$ which implies that $\overline{v}$ is an epimorphism in $\C/\X$.

Note that $\overline{w}=0$ since $\overline{v}$ is an epimorphism.
Thus there exists morphisms $s'\colon \tau^{-1}M\to X_2$ and $t'\colon X_2\to M[1]$ such that
$w=t's'$ where $X_2\in\X$. Note that $t'\in\Hom_{\C}(X_2,M[1])=0$, we have $w=0$.
Therefore the Auslander-Reiten triangle (\ref{t2}) splits, a contradiction.  So we get $M\in\X$.

This shows that  $\{M\in\C~|~\Hom_{\C}(\X,M[1])=0\}\subseteq\X$ and hence $\X$ is cluster tilting.  \qed
\medskip

This theorem immediately yields the following important conclusion.

\begin{cor}
 Let $\C$ be a $2$-Calabi-Yau triangulated category and $\X$ a
non-zero contravariantly finite rigid subcategory of $\C$. Then $\X$ is cluster tilting
if and only if the quotient category $\C/\X$ is abelian.
\end{cor}

Panyue Zhou\\
College of Mathematics, Hunan Institute of Science and Technology, 414006 Yueyang, Hunan, People's Republic of China.\\[1mm]
Email: panyuezhou@163.com


\begin{thebibliography}{99}
\bibitem[AR]{AR} M. Auslander, I. Reiten. \newblock Applications of contravariantly finite subcategories.
\newblock Adv. Math. \textbf{86}(1): 111-152, 1991.



\bibitem[Be]{Be}
A. Beligiannis. Rigid objects, triangulated subfactors and abelian localizations. Math. Z. \textbf{274}: 841-883, 2013.


\bibitem[BK]{BK}
A. I. Bondal and M. M. Kapranov.
\newblock Representable functors, Serre functors, and mutations.
\newblock Mathematics of the USSR-Izvestiya \textbf{35}(3), 519-541, 1990.

\bibitem[BM]{BM}  A. Beligiannis, N. Marmaridis. Left triangulated categories arising from contravariantly
finite subcategories. Comm. Algebra \textbf{22}(12): 5021-5036, 1994.

\bibitem[Ha]{Ha}
D. Happel.
\newblock Triangulated categories in the representation theory of finite-dimensional algebras.
\newblock  London Mathematical Society Lecture Note Series, \textbf{119}.  Cambridge University Press, Cambridge, 1988.


\bibitem[J]{J} P. J{\o}rgensen. Quotients of cluster categories. Proc. Roy. Soc. Edinburgh Sect. A \textbf{140}(1): 65-81, 2010.


\bibitem[KZ]{KZ}
S. Koenig, B. Zhu.
\newblock From triangulated categories to abelian categories: cluster tilting in a general framework.
\newblock Math. Z. 258: 143-160, 2008.

\bibitem[RV]{RV} I. Reiten, M. Van den Bergh. Noetherian hereditary abelian categories satisfying Serre duality.
 J. Amer. Math. Soc. \textbf{15}(2): 295-366, 2002.



\end{thebibliography}
\end{document}